\documentclass[a4paper, 11pt]{extarticle}
\usepackage{cite}
\usepackage{setspace}
\usepackage{mathtext}
\usepackage{amscd}
\usepackage[T2A]{fontenc}
\usepackage[english]{babel}
\usepackage{latexsym,amsfonts,amssymb,amsmath,longtable,amsthm}
\usepackage{bbm,makeidx}
\usepackage[centerlast,small]{caption}
\usepackage{amsmath}
\usepackage[cp1251]{inputenc}
\usepackage[usenames, dvipsnames]{color}
\newtheorem{lemma}{{\scshape Lemma}}
\newtheorem{corollary}{{\scshape Corollary}}
\newtheorem{theorem}{{\scshape Theorem}}
\newtheorem{proposition}{{\scshape Proposition}}

\begin{document}
\title{Twisted conjugacy classes in unitriangular groups}
\author{Timur Nasybullov\footnote{KU Leuven KULAK, 
Etienne Sabbelaan 53, 8500 Kortrijk, Belgium, 
 timur.nasybullov@mail.ru}~\footnote{The author is supported by the Research Foundation -- Flanders (FWO), app. 12G0317N.}}
\date{}
\maketitle
\begin{abstract}
Let $R$ be an integral domain of zero characteristic. In this note we study the Reidemeister spectrum of the group ${\rm UT}_n(R)$ of unitriangular matrices over $R$. We prove that if $R^+$ is finitely generated and $n>2|R^*|$, then  ${\rm UT}_n(R)$ possesses the $R_{\infty}$-property, i.~e. the Reidemeister spectrum of ${\rm UT}_n(R)$ contains only $\infty$, however, if $n\leq|R^*|$, then the Reidemeister spectrum of ${\rm UT}_n(R)$ has nonempty intersection with $\mathbb{N}$. If $R$ is a field, then we prove that the Reidemeister spectrum of ${\rm UT}_n(R)$ coincides with $\{1,\infty\}$, i.~e. in this case ${\rm UT}_n(R)$ does not possess the $R_{\infty}$-property.\\

\end{abstract}
\section{Introduction}
Let $G$ be a group and $\varphi$ be an automorphism of $G$. Two elements $x, y$ of $G$ are said to be (twisted) $\varphi$-conjugated if there exists an element $z \in G$ such that
$x = zy\varphi(z)^{-1}$. The relation of $\varphi$-conjugation is an equivalence relation on $G$ and it divides the group into $\varphi$-conjugacy classes. The number $R(\varphi) \in \mathbb{N}\cup\{\infty\}$ of this classes is called the Reidemeister number of the automorphism $\varphi$. 

Twisted conjugacy classes appear naturally in Nielsen-Reidemeister fixed point theory. Let $X$ be a finite polyhedron and $f:X\to X$ be a homeomorphism of $X$. Two fixed points $x,y$ of $f$ are said to belong to the same fixed point class of $f$ if there exists a path $c:[0,1]\to X$ with $c(0)=x$ and $c(1)=y$ such that $c\simeq f\circ c$, where $\simeq$ denotes the homotopy with fixed endpoints. The relation of being in the same fixed point class is an equivalence relation on the set of fixed points of $f$. The number $R(f)$ of fixed point classes of $f$ is called the Reidemeister number of the homeomorphism $f$. The Reidemeister number $R(f)$ is the homotopic invariant of $f$ and it plays a crucial role in the Nielsen-Reidemeister fixed point theory \cite{Kia}. Denote by $G=\pi_1(X)$ the fundamental group of $X$ and by $\varphi$ the automorphism of $G$ induced by $f$. In this notation the number $R(f)$ of fixed point classes of $f$ is equal to the number $R(\varphi)$ of $\varphi$-conjugacy classes in $G$ (see \cite[Chapter III, Lemma 1.2]{Kia}). Thus, the topological problem
of finding $R(f)$ reduces to the purely algebraic problem of finding $R(\varphi)$. If $X$ is a finite polyhedron, then $\pi_1(X)$ is a finitely presented group, therefore twisted conjugacy classes in finitely presented groups are of special interest.

The Reidemeister spectrum $Spec_R(G)$ of a group $G$ is the subset of $\mathbb{N}\cup \{\infty\}$ of the form $Spec_R(G)=\{R(\varphi)~|~\varphi\in{\rm Aut}(G)\}$. If $Spec_R(G) = \{\infty\}$, then the group $G$ is said to possess the $R_{\infty}$-property. The problem of classifying groups
which possess the $R_{\infty}$-property was proposed by A.~Fel'shtyn and R.~Hill in \cite{FelHil}. The study of this problem
has been quite an active research topic in recent years. We refer to the paper \cite{FelNas} for an overview of the families of groups which have been studied in this context until 2016. More recent results can be found in \cite{DekGon2, MubSan, Tro, Tro2}. The author studied twisted conjugacy classes and the $R_{\infty}$-property for classical linear groups \cite{FelNas, Nas, Nas2, Nas3, Nas4}. For the immediate consequences of the $R_{\infty}$-property for topological fixed point
theory see \cite{GonWon}. Some aspects of the $R_{\infty}$-property can be found in \cite{FelTro2}.

The Reidemeister numbers for automorphisms of nilpotent groups are studied only for free nilpotent groups \cite{DekGon, Rom}, nilpotent quotients of surface groups \cite{DekGon2} and some very specific nilpotent groups motivated by geometry \cite{GonWon}. Denote by $N_{r,c}=F_r/\gamma_{c+1}(F_r)$ the free nilpotent group of rank $r$ and nilpotency class $c$. V.~Roman'kov proved in \cite{Rom} that if $r\geq 4$ and $c\geq 2r$, then the group $N_{r,c}$ possesses the $R_{\infty}$-property. He also found the Reidemeister spectrum for groups $N_{2,2}$, $N_{2,3}$, $N_{3,2}$. K.~Dekimpe and D.~Gon\c{c}alves extended the result of Roman'kov in \cite{DekGon} proving that the group $N_{r,c}$ possesses the $R_{\infty}$-property if and only if $r\geq 2$ and $c\geq 2r$. The authors of \cite{DekTerVar} studied the Reidemeister spectrum for the groups $N_{r,c}$ for $c<2r$. In this note we study the Reidemeister spectrum and the $R_{\infty}$-property for another class of  nilpotent groups, namely for groups ${\rm UT}_n(R)$ of upper unitriangular matrices over integral domains $R$ of zero characteristic.  The group ${\rm UT}_n(\mathbb{Z})$ is very important for studying nilpotent groups in general since if $G$ is a finitely generated torsion free nilpotent group, then it can be imbedded into ${\rm UT}_n(\mathbb{Z})$ for an appropriate positive integer $n$. The group ${\rm UT}_n(R)$ is both nilpotent and classical linear, so, the present work is a natural continuation of works \cite{FelNas, Nas, Nas2, Nas3, Nas4}. The main results of the paper are formulated in the following two theorems.

~\\
\noindent	\textbf{{\scshape Theorem 1.}} Let $R$ be an infinite integral domain of zero characteristic with finitely generated additive group. If $n > 2|R^*|$, then ${\rm UT}_n(R)$ possesses the $R_{\infty}$-property.

~\\
\noindent	\textbf{{\scshape Theorem 2.}} Let $R$ be a field of zero characteristic and $n \geq 3$. Then the Reidemeister spectrum of ${\rm UT}_n(R)$ is $\{1,\infty\}$.

~\\
The author is grateful to V. Roman'kov for the attention to the work and important remarks about the automorphism group of the group of unitriangular matrices.
\section{Preliminaries}
In this section we will introduce preliminary facts which we need for proving the main results of the paper. 
The following proposition can be found, for example, in \cite[Corollary 2.5]{FelTro}.
\begin{proposition}\label{inn} Let $G$ be a group, $\varphi$ be an automorphism and $\psi$ be an inner
automorphism of $G$. Then $R(\psi\varphi) = R(\varphi)$.
\end{proposition}
The following two propositions are obvious.
\begin{proposition}\label{ind}Let $G$ be an abelian group and $\varphi$ be an automorphism of $G$. Then
the $\varphi$-conjugacy class $[e]_{\varphi}$ of the unit element $e$ is a subgroup of $G$ and $R(\varphi) = \big|G : [e]_{\varphi}\big|$.
\end{proposition}
\begin{proposition}\label{prod}Let $G$ be a direct product of groups $P$ and $Q$ and $\varphi$ be an automorphism of $G$ such that $\varphi(P) = P$, $\varphi(Q) = Q$. Then $R(\varphi) = R(\varphi|_P)R(\varphi|_Q)$.
\end{proposition}
The following statement is proved in \cite[Lemma 1.1(3)]{GonWon}.
\begin{proposition}\label{zf}
Let $G$ be a group, $\varphi$ be an automorphism of $G$ and $H$ be a $\varphi$-admissible central subgroup of $G$. Denote by $\varphi^{\prime}, \overline{\varphi}$ respectively the automorphism of $H$ and the automorphism of $G/H$ induced by $\varphi$. Then $R(\varphi)\leq R(\varphi^{\prime})R(\overline{\varphi})$.
\end{proposition}
The original statement \cite[Lemma 1.1 (3)]{GonWon} says that in the conditions of Proposition \ref{zf} there is an equality $R(\varphi) = R(\varphi^{\prime})R(\overline{\varphi})$. However the proof from \cite[Lemma 1.1 (3)]{GonWon} guarantees only the inequality. Moreover, the statement with equality $R(\varphi) = R(\varphi^{\prime})R(\overline{\varphi})$ is not correct. For example, if $G=\langle x~|~x^4=1\rangle$ is a cyclic group of order $4$, $H=\langle x^2\rangle$ is a subgroup of $G$ and $\varphi:x\mapsto x^{-1}$ is an automorphism of $G$, then $R(\varphi)=R(\varphi^{\prime})=R(\overline{\varphi})=2$ and $R(\varphi) \neq R(\varphi^{\prime})R(\overline{\varphi})$.
Proposition \ref{zf} has the following corollary.
\begin{corollary}\label{nilin} Let $G$ be a nilpotent group of class $n$ and let $1 = Z_0 < Z_1 < \dots <
Z_n = G$ be the upper central series of $G$. For an automorphism $\varphi$ of $G$ denote by
$\varphi_k$ the automorphism of $Z_{k+1}/Z_k$ induced by $\varphi$. Then $R(\varphi) \leq \prod_{k=1}^{n-1} R(\varphi_k)$.
\end{corollary}
If $G$ is a finitely generated torsion free nilpotent group, then Corollary \ref{nilin} has the following stronger formulation proved in \cite[Lemma 2.7]{Rom}.
\begin{proposition}\label{mrom}Let $G$ be a finitely generated torsion free nilpotent group of class
$n$ and let $1 = Z_0 < Z_1 < \dots < Z_n = G$ be the upper central series of $G$. For an automorphism $\varphi$ of $G$ denote by $\varphi_k$ the automorphism of $Z_{k+1}/Z_k$ induced by $\varphi$.
Then the number $R(\varphi)$ is finite if and only if the numbers $R(\varphi_k)$ are finite for all $k = 0,\dots, n-1$. Moreover, in this case $R(\varphi)=\prod_{k=0}^{n-1}R(\varphi_k)$.
\end{proposition}
The following statement can be found in \cite[Lemma 4.1]{Rom}.
\begin{proposition}\label{fixab}Let $\varphi$ be an automorhism of a free abelian group of finite rank.
Then $R(\varphi)$ is infinite if and only if $\varphi$ has a non-trivial fixed point.
\end{proposition}
An associative, commutative ring $R$ with unit and without zero divisors is called
an integral domain. We denote by $R^+$ the additive group of $R$, and by $R^*$ the
multiplicative group of $R$. Proposition \ref{fixab} has the following corollary.
\begin{corollary}\label{inR}
Let $R$ be an infinite integral domain such that $R^+$ is a finitely
generated abelian group. Let $a\neq 1$ be an element from $R^{*}$ and $\varphi: x\mapsto ax$ be an
automorphism of $R^+$. Then $R(\varphi)$ is finite.
\end{corollary}

Let $R$ be an integral domain of zero characteristic. For $n \geq 2$ denote by ${\rm M}_n(R)$ the set of all $n \times n$ matrices with entries from $R$, and by ${\rm I}_n$ the identity matrix in
${\rm M}_n(R)$. For $1 \leq i, j \leq n$ denote by ${\rm E}_{i,j}$ the matrix which has $1$ in the position
$(i, j)$ and $0$ in all other positions, and for $1 \leq i <j \leq n$,  $x \in R$ denote by ${\rm T}_{i,j}(x) = {\rm I}_n + x{\rm E}_{i,j}$. The group ${\rm UT}_n(R)$ of upper unitriangular matrices over $R$ is the set
of matrices from ${\rm M}_n(R)$ which have $1$ on the diagonal, $0$ below the diagonal, and
arbitrary entries above the diagonal
$${\rm UT}_n(R)=\Big\{(x_{i,j})\in{\rm M}_n(R)~\Big|~x_{i,i}=1 ~\text{for}~ 1\leq i\leq n,~x_{i,j}=0 ~\text{for}~ 1\leq j<i\leq n\Big\}.$$
This group is generated by the elements ${\rm T}_{i,i+1}(x)$ for $i = 1, \dots, n-1$, $x \in R$.

The group ${\rm UT}_n(R)$ is nilpotent of class $n-1$. The $k$-th member $Z_k({\rm UT}_n(R))$ of the upper central series of ${\rm UT}_n(R)$ consists of the matrices which have zeros on the
first $n-1-k$ superdiagonals
$$Z_k=Z_k({\rm UT}_n(R))=\Big\{(x_{i,j})\in{\rm UT}_n(R)~\Big|~x_{i,j}=0 ~\text{for}~ 0<j-i<n-k\Big\}.$$
The quotient $Z_{k+1}/Z_k$ is isomorphic to the direct sum
$$Z_{k+1}/Z_k=(R^{+})^{k+1}=\{x_1e_1+\dots+x_{k+1}e_{k+1}~|~x_1,\dots,x_{k+1}\in R^{+}\},$$
where $e_r$ denotes the column of length $k + 1$ which has $1$ in the position $r$ and $0$ in all other positions.

We recall the definitions of some classes of automorphisms of the group ${\rm UT}_n(R)$ over an integral domain $R$ (see \cite{Lev, Mah}).

For a matrix $A \in {\rm UT}_n(R)$ denote by $\varphi_A$ the automorphism of ${\rm UT}_n(R)$ acting by
the rule $\varphi_A$(X) = $AXA^{-1}$. The automorphism $\varphi_A$ is called an inner automorphism
of ${\rm UT}_n(R)$. All inner automorphisms generate a normal subgroup in the group of all
automorphisms of ${\rm UT}_n(R)$.

For a diagonal matrix $D = {\rm diag}(d_1,\dots,d_n)$, where $d_1,\dots,d_n \in R^*$ denote by  $\psi_D$ the automorphism of ${\rm UT}_n(R)$ acting by the rule  $\psi_D(X)$ = $DXD^{-1}$. The automorphism
$\psi_D$ is called the diagonal automorphism. It induces the automorphism $\psi_{D,k}$ of $Z_{k+1}/Z_k$ which acts by the rule  $\psi_{D,k}(x_re_r) = d_rd_{n+r-(k+1)}^{-1}x_re_r$ for $r \leq k + 1$.

For a homomorphism $\lambda : R^+\to R^+$ denote by $\Lambda$ the automorphism of ${\rm UT}_n(R)$ which acts on generators of ${\rm UT}_n(R)$ by the rule $\Lambda({\rm T}_{i,i+1}(x)) = {\rm T}_{i,i+1}(x){\rm T}_{1,n}(\lambda(x))$. This automorphism acts identically modulo the center of ${\rm UT}_n(R)$ and is called the central automorphism of ${\rm UT}_n(R)$. If $n \geq 3$, then every central automorphism acts trivially on all quotients $Z_{k+1}/Z_k$.

For an automorphism $\delta$ of the ring $R$ denote by $\Delta$ the map from ${\rm M}_n(R)$ to itself which acts on the matrix $X = (x_{i,j})$ by the rule $\Delta(X) = (\delta(x_{i,j}))$. Delta induces
the automorphism of ${\rm UT}_n(R)$ which is called the ring automorphism. If  $\psi_D$ is a
diagonal automorphis and $\Delta$ is a ring automorphism of ${\rm UT}_n(R)$, then we have the equality $\Delta\psi_D=\psi_{\Delta(D)}\Delta$.

Denote by $\sigma$ the map from ${\rm M}_n(R)$ to itself which maps the matrix $(x_{i,j})$ to the matrix $(x_{n-j+1,n-i+1})^{-1}$. This  map is given by flipping the matrix by the antidiagonal and then taking the inverse  induces the antiautomorphism of ${\rm UT}_n(R)$ which is called the flip automorphism. The order of $\sigma$ is equal to $2$ and $\sigma$ commutes with every ring automorphism. If $\psi_D$ is a diagonal automorphis of ${\rm UT}_n(R)$, then we have the
equality $\sigma\psi_D\sigma=\psi_{\sigma(D)}$.

The following result which gives the discription of an arbitrary automorphism of ${\rm UT}_n(R)$ is a consequence of \cite[Theorems 1, 2, 3]{Lev} in assumption that $R$ is an integral domain of zero characteristic. 
\begin{proposition}\label{lev}Let $R$ be an integral domain, $n\geq 3$ be a positive integer and $\varphi$ be an automorphism of ${\rm UT}_n(R)$, then there exist an inner automorphism $\varphi_A$, a diagonal automorphism
$\psi_D$, a central automorphism $\Lambda$, a ring automorphism $\Delta$ and an integer $m\in\{0,1\}$
such that $\varphi=\varphi_A\Lambda\sigma^m\psi_D\Delta$.
\end{proposition}
\section{$R$ is a ring}
In this section we will study the $R_{\infty}$-property for the group ${\rm UT}_n(R)$ in the case when $R$ is an infinite integral domain such that $R^+$ is a finitely generated abelian group. In this case ${\rm UT}_n(R)$ is finitely presented. From now on we denote by $Z_k = Z_k({\rm UT}_n(R))$. In order to denote the $k$-th member of the upper central series of some group $G\neq {\rm UT}_n(R)$ we will write $Z_k(G)$.
\begin{proposition}\label{notinf}Let $R$ be an infinite integral domain such that $R^+$ is a finitely
generated abelian group. If $n \leq |R^*|$, then the group ${\rm UT}_n(R)$ does not possess the $R_{\infty}$-property.
\end{proposition}
\noindent \textbf{Proof.} Since $n \leq |R^*|$, there exist $n$ different invertible elements $d_1,\dots,d_n$ from $R$. Denote by $D = {\rm diag}(d_1,\dots, d_n)$, $\varphi=\psi_D$ and let us prove that $R(\varphi)$ is finite.

By Proposition \ref{mrom} it is enough to prove that for all $k$ the numbers $R(\varphi_k)$ are
finite, where $\varphi_k$ denotes the automorphism of $Z_{k+1}/Z_k = (R^+)^{k+1} = \{e_1x_1 +\dots +
e_{k+1}x_{k+1}~|~x_1,\dots,x_{k+1}\in R^+\}$ induced by $\varphi$. The automorphism $\varphi_k$ acts by the rule
$\varphi_k(x_re_r) = d_rd_{n-(k+1)+r}^{-1}x_re_r$ for $r \leq k+1$. Since $\varphi_k(R^+e_r) = R^+e_r$, by Proposition~\ref{prod}
it is enought to prove that all numbers $R(\varphi_{k,r})~(r = 1,\dots, k + 1)$ are finite, where
$\varphi_{k,r}$ is the automorphism of $R^+e_r$ induced by $\varphi_k$. Since all the elements $d_1,\dots,d_n$
are different, all the elements $d_rd_{n-(k+1)+r}^{-1}$ are not equal to $1$. The automorphism $\varphi_{k,r}$ acts by the rule $\varphi_{k,r}(x_re_r) = d_rd_{n-(k+1)+r}^{-1}x_re_r$, therefore by Corollary~\ref{inR} the
number $R(\varphi_{k,r})$ is finite. Therefore the numbers $R(\varphi_k)$ are finite and the number
$R(\varphi$) is finite.\hfill$\square$

If $p$ is a positive integer which is not a perfect square, then the ring $\mathbb{Z}[\sqrt{p}]$ has infinitely many invertible elements. Therefore by Proposition \ref{notinf} the group ${\rm UT}_n(\mathbb{Z}[\sqrt{p}])$ does not possess the $R_{\infty}$-property independently on $n$.
\begin{theorem}\label{th1}Let $R$ be an infinite integral domain such that $R^+$ is a finitely generated abelian group. If $n > 2|R^*|$, then ${\rm UT}_n(R)$ possesses the $R_{\infty}$-property.
\end{theorem}
\noindent\textbf{Proof.} Let $\varphi$ be an automorphism of ${\rm UT}_n(R)$. We are going to prove that the
number $R(\varphi)$ is infinite. By Proposition \ref{lev} there exist an inner automorphism $\varphi_A$, a diagonal automorphism $\psi_D$, a central automorphism $\Lambda$, a ring automorphism $\Delta$
and an integer $m \in\{0,1\}$ such that $\varphi=\varphi_A\Lambda\sigma^m\psi_D\Delta$. By Proposition \ref{inn} we can
assume that $\varphi_A = id$, so, $\varphi=\Lambda\sigma^m\psi_D\Delta$. Since $R$ is an infinite integral domain and
$R^+$ is a finitely generated abelain group, the characteristic of $R$ is equal to zero, therefore $R^*$ contains at least two elements $(\pm1)$ and $n > 4$. By Proposition~\ref{mrom} it
is enough to prove that there exists a number $k$ such that $R(\varphi_k) = \infty$, where $\varphi_k$
is the automorphism of $Z_{k+1}/Z_k$ induced by $\varphi$. Since $n > 4$, the automorphism $\Lambda$
acts trivially on all quotients $Z_{k+1}/Z_k$ and therefore we can assume that $\Lambda=id$ and
$\varphi=\sigma^m\psi_D\Delta$. Let $D = {\rm diag}(d_1,\dots, d_n)$ be a diagonal matrix. Depending on the number $m \in\{0,1\}$ we have two cases.

\textit{Case 1: $\varphi=\psi_D\Delta$.} Since $n > |R^*|$, there exist two indices $1 \leq j < i \leq n$ such
that $d_i = d_j$. Denote by $k + 1 = n + j - i$ and let us prove that $R(\varphi_k) = \infty$, where
$\varphi_k$ is an automorphism of $Z_{k+1}/Z_k = \{e_1x_1 + \dots + e_{k+1}x_{k+1}~|~ x_1,\dots,x_{k+1}\in R^+\}$ induced by $\varphi$.

The automorphism $\varphi_k$ acts by the rule $\varphi_k : x_re_r \mapsto d_rd_{n-(k+1)+r}^{-1}\delta(x_r)e_r$. Since
$\varphi_k(R^+e_r) = R^+e_r$, by Proposition \ref{prod} it is sufficient to prove that $R(\varphi_{k,r}) = \infty$ for
some $r \leq k+1$, where $\varphi_{k,r}$ denotes the automorphism of $R^+e_r$ induced by $\varphi_k$. Since
$k +1 = n+j -i$, we have $j = k +1-(n-i) \leq k +1$ and we can speak about $\varphi_{k,j}$. This automorphism acts on $R^+e_j$ by the rule
$$\varphi_{k,j}(xe_j)=d_jd_{n-(k+1)+j}^{-1}\delta(x)e_j=d_jd_i^{-1}\delta(x)e_j=\delta(x)e_j.$$
This automorphism obviously has a non-trivial fixed point ($1e_j \neq 0e_j$), therefore by Proposition~\ref{fixab} we have $R(\varphi_{k,j}) = \infty$, then $R(\varphi_k) = \infty$ and $R(\varphi) = \infty$.

\textit{Case 2: $\varphi=\sigma\psi_D\Delta$.} Denote by $N$ the rank of the finitely generated abelian group $R^+$ and consider an arbitrary element $x \in R$. There exist integers $\alpha_0,\dots,\alpha_N$ not all equal to zero such that $\alpha_0+\alpha_1x+\dots+\alpha_Nx^N=0$. Since $\delta$ is an automorphism
of $R$ and $\alpha_0,\dots,\alpha_N$ are integers, we have $\alpha_0+\alpha_1\delta(x)+\dots+\alpha_N\delta(x)^N=0$. Therefore
all the elements $x, \delta(x), \delta^2(x),\dots$ are roots of the polynomial $\alpha_0+\alpha_1T +\dots+\alpha_NT^N$
and $\delta$ has a finite order which divides $(N + 1)!$.

Denote by $A = \sigma(D)\Delta(D) = {\rm diag}(a_1,\dots,a_n)$, where $a_r = \delta(d_r)d_{n+1-r}^{-1}$ for
$r = 1,\dots,n$. Since $n > 2|R^*|$ there exist $3$ indices $i, j, l$ such that $a_i = a_j = a_l$. Let $1 \leq j < i \leq n$ be two of this indices which satisfy the condition $i + j \neq n + 1$.
Denote by $k + 1 = n + j - i$ and let us prove that $R(\varphi_k) = \infty$, where $\varphi_k$ is the
automorphism of $Z_{k+1}/Z_k = \{x_1e_1 +\dots+x_{k+1}e_{k+1}~|~x_1,\dots,x_{k+1} \in R^+\}$ induced by $\varphi$. The automorphism $\varphi_k$ acts by the rule
\begin{equation}\label{f1}
\varphi_k: x_re_r \mapsto -d_rd_{n-(k+1)+r}^{-1}\delta(x_r)e_{k+2-r}
\end{equation}
for $r = 1,\dots,k+1$. Since $i+j \neq n+1$, we have $k+2-j = n+j-i+1-j = n+1-i \neq j$. Consider the subgroup $H = R^+e_j \oplus R^+e_{k+2-j}$ of $Z_{k+1}/Z_k$. The automorphism $\varphi_k$
fixes $H$, so, denote by $\theta$ the automorphism of $H$ induced by $\varphi_k$. By Proposition \ref{prod} it is enough to prove that $R(\theta) = \infty$.

Since the characteristic of $R$ is equal to zero, $R$ contains the set of integers. Suppose that the number of $\theta$-conjugacy classes in $H$ is finite and consider the set
of elements $y_t = 2te_j + t\delta^{-1}(d_{k+2-j}^{-1}d_{n+1-j})e_{k+2-j}$ for $t = 0,1,2,\dots$ Since $R(\theta)$ is finite
there are two numbers $r, s$ such that $y_r$ and $y_s$ are $\theta$-conjugate, i. e. there exists an element $x \in H$ such that $y_r = x + y_s -\theta(x)$, or equivalently, $y_t = x -\theta(x)$, where
$t = r -s\neq 0$. From this equality follows that
\begin{equation}\label{f2}
y_t + \theta(y_t) = x -\theta^2(x).
\end{equation}
Note that since $\delta$ is an automorphism of $R$ and $t$ is an integer, we have $\delta(t) = t$. Let us rewrite $y_t + \theta(y_t)$ in more details using equality (\ref{f1})
\begin{align}
\notag y_t+\theta(y_t)&=2te_j+t\delta^{-1}(d_{k+2-j}^{-1}d_{n+1-j})e_{k+2-j}+\theta\Big(2te_j+t\delta^{-1}(d_{k+2-j}^{-1}d_{n+1-j})e_{k+2-j}\Big)\\
\notag&=2te_j+t\delta^{-1}(d_{k+2-j}^{-1}d_{n+1-j})e_{k+2-j}\\	
\notag&~~-2d_jd_{n-(k+1)+j}^{-1}te_{k+2-j}-d_{k+2-j}d_{n+1-j}^{-1}\delta\Big(t\delta^{-1}(d_{k+2-j}^{-1}d_{n+1-j})\Big)e_j\\
\label{f3}&=te_j+ce_{k+2-j},
\end{align}
where $c=t\Big(\delta^{-1}(d_{k+2-j}^{-1}d_{n+1-j})-2d_jd_{n-(k+1)+j}^{-1}\Big)$. Since $\varphi=\sigma\psi_D\Delta$, we have
\begin{align}
\notag\varphi^2&=\sigma\psi_D\Delta\sigma\psi_D\Delta=\sigma\psi_D\sigma\Delta\psi_D\Delta=\psi_{\sigma(D)}\Delta\psi_D\Delta\\
\label{f4}&=\psi_{\sigma(D)}\psi_{\Delta(D)}\Delta^2=\psi_{\sigma(D)\Delta(D)}\Delta^2=\psi_A\Delta^2.
\end{align}
Let the element $x$ from equality (\ref{f2}) has the form $x = x_1e_j+x_2e_{k+2-j}$ for $x_1, x_2 \in R^+$.
Therefore by equality (\ref{f4}) the element $\theta^2(x)$ has the form
\begin{align}
\notag\theta^2(x)&=a_ja_{n+j-(k+1)}^{-1}\delta^2(x_1)e_j + a_{k+2-j}a_{n+k+2-j-(k+1)}^{-1}\delta^2(x_2)e_{k+2-j}\\
\notag&=a_ja_i^{-1}\delta^2(x_1)e_j + a_{k+2-j}a_{n+1-j}^{-1}\delta^2(x_2)e_{k+2-j}\\
\label{f5}&=\delta^2(x_1)e_j + a_{k+2-j}a_{n+1-j}^{-1}\delta^2(x_2)e_{k+2-j}.
\end{align}
Using equalities (\ref{f3}) and (\ref{f5}) we can rewrite equality (\ref{f2}) in more details
$$\delta^2(x_1)e_j + a_{k+2-j}a_{n+1-j}^{-1}\delta^2(x_2)e_{k+2-j} = \theta^2(x) = (x_1 - t)e_j + (x_2 - c)e_{k+2-j}.$$
This equality implies that $\delta^2(x_1) = x_1 - t$ and $\delta^{2p}(x_1) = x_1 - pt$ for an arbitrary
positive integer $p$. This contradicts the fact that $\delta$ is an automorphism of $R$ of finite
order. Therefore there are infinitely many $\theta$-conjugacy classes in $H$, $R(\varphi_k) = \infty$ and $R(\varphi) = \infty$.\hfill$\square$

If $p \neq 1$ is a non-negative integer, then the ring $\mathbb{Z}[i\sqrt{p}]$ has only two invertible elements $(\pm 1)$, therefore by Theorem \ref{th1} for $n \geq 5$ the group ${\rm UT}_n(\mathbb{Z}[i\sqrt{p}])$
possesses the $R_{\infty}$-property. The ring $\mathbb{Z}[i]$ of Gauss integers has 4 invertible elements $(\pm1,\pm i)$. Therefore by Theorem \ref{th1} for $n \geq 9$ the group ${\rm UT}_n(\mathbb{Z}[i])$ possesses the $R_{\infty}$-property.
However, by Proposition \ref{notinf} the groups ${\rm UT}_n(\mathbb{Z}[i])~(n = 2, 3, 4)$ do not possess the $R_{\infty}$-property.

Proposition \ref{notinf} and Theorem \ref{th1} do not cover the situation when $|R^{*}|<n\leq 2|R^*|$. The group ${\rm UT}_3(\mathbb{Z})$ is isomorphic to the free nilpotent group $N_{2,2}$ of rank $2$ and nilpotency class $2$, which has the Reidemeister spectrum $Spec_R(N_{2,2})=2\mathbb{N}\cup\{\infty\}$ (see \cite{Rom}), i.~e. ${\rm UT}_3(\mathbb{Z})$ does not possess the $R_{\infty}$-property. The group ${\rm UT}_4(\mathbb{Z})$ is known to possess the $R_{\infty}$-property \cite[Theorem 9.2.11]{Dug}. Since $\mathbb{Z}$ has only $2$ invertible elements, both $n=3$ and $n=4$ satisfy the condition $|\mathbb{Z}^{*}|<n\leq 2|\mathbb{Z}^*|$, therefore for $|R^{*}|<n\leq 2|R^*|$ both situations are possible: the group ${\rm UT}_n(R)$ can either possess or not possess the $R_{\infty}$-property.
\section{R is a field}
In this section we will study the Reidemeister spectrum of the group ${\rm UT}_n(R)$ in the
case when $R$ is a field of zero characteristic.

If $R$ is a field and $a \neq 1$ is an element from $R$, then the map $\varphi: x \mapsto ax$ is the automorphism of $R^+$ with $R(\varphi) = 1$. Using this remark and Corollary \ref{nilin} (instead of Proposition \ref{mrom}) by the same arguments as in the proof of Proposition \ref{notinf} we conclude
that there exists a diagonal automorphism $\varphi$ of ${\rm UT}_n(R)$ with $R(\varphi) = 1$. We are
going to prove that $1$ is the only possible finite value of $R(\varphi)$ in ${\rm UT}_n(R)$.
\begin{lemma}\label{lem1} Let $R$ be an integral domain of zero characteristic, $n \geq 3$ and $\varphi=\Lambda\sigma^m\psi_D\Delta$ be an
automorphism of ${\rm UT}_n(R)$. If matrices $A,B$ from $Z_{k+1}$
are $\varphi$-conjugated, then there exists a matrix $Y \in Z_{k+1}$ and a matrix $Z \in {\rm UT}_n(R)$ such that $\varphi(Z) = \Lambda(Z)$ and $A = (Y Z)^{-1}B\varphi(Y Z)$.
\end{lemma}
\noindent \textbf{Proof.} Denote by $\theta_r$ the natural homomorphism ${\rm UT}_n(R) \to {\rm UT}_n(R)/Z_r$ and by $\overline{\varphi}_r$ the automorphism of ${\rm UT}_n(R)/Z_r$ induced by $\varphi$. Since $A,B$ are $\varphi$-conjugated, there
exists an element $X \in {\rm UT}_n(R)$ such that
\begin{equation}\label{f6}
A=X^{-1}B\varphi(X).
\end{equation}
The element $X$ belongs to $Z_r$ for some $r \geq k+1$. In order to prove the statement we
will use induction on the number $r-(k+1)$. The basis is obvious: if $r-(k+1) = 0$, then $X \in Z_{k+1}$, so, denote by $Y = X$, $Z = {\rm I}_n$.

If $X \in Z_{r+1}$, then $X$ can be expressed in the form $X = PQ$, where $P \in Z_r$ and $Q$ is a typical representative of the element from $Z_{r+1}/Z_r$ in ${\rm UT}_n(R)$, i. e. the matrix which differs
from ${\rm I}_n$ only in the $(n-r)$-th superdiagonal. Acting on equality (\ref{f6}) by the homomorphism
$\theta_r$ we have $\overline{\varphi}_r(\theta_r(Q)) = \theta_r(Q)$, or in details $\sigma^m\psi_D\Delta(\theta_r(Q)) = \theta_r(Q)$ ($\Lambda$ disappeared
since $\Lambda$ acts identically modulo $Z_1$ and we factorized by $Z_r$). Since $Q$ is a typical
representative of the element from $Z_{r+1}/Z_r$, we conclude that $\sigma^m\psi_D\Delta(Q)=Q$ or $\varphi(Q) = \Lambda(Q)$. So, from equality (\ref{f6}) we have $A = Q^{-1}(P^{-1}B\varphi(P))\varphi(Q)$ or
$$QA\varphi(Q)^{-1}=P^{-1}B\varphi(P),$$
for $P \in Z_m$ and $Q$ satisfying $\varphi(Q) = \Lambda(Q) = QC$ for some $C \in Z_1$. The matrix
$QA\varphi(Q)^{-1} = QAQ^{-1}C^{-1}$ belongs to $Z_{k+1}$ and the matrix $P$ belongs to $Z_m$, therefore
we can use the induction conjecture, i. e. there exists a matrix $Y^{\prime} \in Z_{k+1}$ and a matrix $Z^{\prime}$ with $\varphi(Z^{\prime}) = \Lambda(Z^{\prime})$ such that $QA\varphi(Q)^{-1} = (Y^{\prime}Z^{\prime})^{-1}B\varphi(Y^{\prime} Z^{\prime}$) or
$$A = (Y^{\prime} Z^{\prime}Q)^{-1}B\varphi(Y^{\prime}Z^{\prime}Q).$$
Denoting by $Y = Y^{\prime}$, $Z = Z^{\prime}Q$ we finish the proof.\hfill $\square$

Lemma \ref{lem1} has two useful corollaries.
\begin{corollary}\label{cor3}Let $R$ be a field of zero characteristic, $n \geq 3$ and $\varphi=\Lambda\sigma^m\psi_D\Delta$ be an automorphism of ${\rm UT}_n(R)$ with $R(\varphi) < \infty$. Then all central elements from
${\rm UT}_n(R)$ are $\varphi$-conjugated.
\end{corollary}
\noindent \textbf{Proof.} Let $D$ be a diagonal matrix of the form $D = {\rm diag}(d_1,\dots,d_n)$. Denote by $H$ the following subset of $R$
$$H = \{\lambda(t_1) + \lambda(t_2) + \dots + \lambda(t_{n-1}) +
\big(y-d_1d_n^{-1}\delta(y)\big)~|~t_i,y\in R,~d_id_{i+1}^{-1}\delta(t_i)=t_i\}$$
if $m=0$ and 
$$H = \{\lambda(t_1) + \lambda(t_2) + \dots + \lambda(t_{n-1}) +
\big(y-d_1d_n^{-1}\delta(y)\big)~|~t_i,y\in R,~d_id_{i+1}^{-1}\delta(t_{n-i})=t_i\}$$
if $m=1$. We are going to prove that the central matrices $A = {\rm T}_{1,n}(a)$ and $B ={\rm T}_{1,n}(b)$ are $\varphi$-conjugated if and only if $a - b \in H$. We consider in details only the
case $m = 0$, i. e. $\varphi=\Lambda\psi_D\Delta$. The case $m = 1$ is similar.

If $a-b \in H$, then for elements $t_1,\dots,t_{n-1},y\in R$ with conditions $d_id_{i+1}^{-1}= t_i$ for $i = 1,\dots,n-1$ we have $a - b = \lambda(t_1) + \lambda(t_2) + \dots + \lambda(t_{n-1}) +
\big(y - d_1d_{n}^{-1} \delta(y)\big)$. Denote by $T = {\rm I}_n + {\rm E}_{1,2}(t_1) + {\rm E}_{2,3}(t_2) + \dots + {\rm E}_{n-1,n}(t_{n-1})$ and by $Y = {\rm T}_{1,n}(y)$. In
this denotation we have $A = (TY)^{-1}B\varphi(TY)$, i. e. $A$ and $B$ are $\varphi$-conjugated.

Conversely, let $A = {\rm T}_{1,n}(a)$ and $B = {\rm T}_{1,n}(b)$ be two central $\varphi$-conjugated matrices from ${\rm UT}_n(R)$. By Lemma \ref{lem1} there exists a matrix $Y \in Z_1$ and a matrix
$Z \in {\rm UT}_n(R)$ such that $\varphi(Z) = \Lambda(Z)$ and $A = (Y Z)^{-1}B\varphi(Y Z)$, or equivalently
\begin{equation}\label{f7}
ZA\varphi(Z)^{-1}=Y^{-1}B\varphi(Y).
\end{equation}
The matrix $Z$ can be expressed as a product $Z = TS$, where $S \in Z_{n-2}$ and $T$ is a
typical representative of $Z_{n-1}/Z_{n-2}$ in ${\rm UT}_n(R)$, i. e. the matrix which differs from ${\rm I}_n$ only in the first superdiagonal. Since $\Lambda$ acts identically on $Z_{n-2}$, we have $\Lambda(S) = S$ and equality (\ref{f7}) can be rewritten in the form
\begin{equation}\label{f8}
TA\Lambda(T)^{-1}=Y^{-1}B\varphi(Y).
\end{equation}
Since $Y\in Z_1$, it has the form $Y={\rm T}_{1,n}(y)$. Since $T$ is a
typical representative of $Z_{n-1}/Z_{n-2}$, it has the form $T = {\rm I}_n + {\rm E}_{1,2}(t_1) + {\rm E}_{2,3}(t_2) + \dots + {\rm E}_{n-1,n}(t_{n-1})$. Therefore from equality (\ref{f8}) we have
$$a-\lambda(t_1)-\lambda(t_2)-\dots-\lambda(t_{n-1})=b+\big(y-d_1d_n^{-1}\delta(y)\big)$$
where elements $t_1,\dots,t_{n-1}$ satisfy the conditions $d_id_{i+1}^{-1}= t_i$ (which follow from
the equality $\varphi(T) = \Lambda(T)$), therefore $a - b \in H$.

The set $H$ is obviously a subgroup of $R^{+}$. The fact that ${\rm T}_{1,n}(a)$ and ${\rm T}_{1,n}(b)$ are $\varphi$-conjugated if and only if $a-b\in H$ implies that the number of non-$\varphi$-conjugated central elements of ${\rm UT}_n(R)$ is equal to the index $\big|R^{+}:H\big|$. Since $R(\varphi)$ is finite, the index $\big|R^{+}:H\big|$ is also finite. Since $R$ is a field of zero characteristic, the group $R^{+}$ is a divisible abelian group, and it has only one subgroup of finite index (itself), i.~e. $H=R^{+}$. Therefore all central elements are twisted $\varphi$-conjugated.\hfill$\square$
\begin{corollary}\label{cor4}Let $R$ be a field of zero characteristic, $n \geq 3$ and $\varphi = \Lambda\sigma^m\psi_D\Delta$ be an automorphism of ${\rm UT}_n(R)$ with $R(\varphi) < \infty$. Then every element from ${\rm UT}_n(R)$ is $\varphi$-conjugated with some element from $Z_1$.
\end{corollary}
\noindent \textbf{Proof.} Denote by $\varphi_k$ the automorphism of $Z_{k+1}/Z_k$ induced by $\varphi$. Denote also
by $G = {\rm UT}_n(R)/Z_1$, by $\psi$ the automorphism of $G$ induced by $\varphi$ and by $\psi_k$ the automorphism of $Z_{k+1}(G)/Z_k(G)$ induced by  $\psi$.

If two elements $A,B \in Z_{k+1}~(k \neq 0)$ are $\varphi$-conjugated, then by Lemma \ref{lem1} there
exists a matrix $Y \in Z_{k+1}$ and a matrix $Z \in {\rm UT}_n(R)$ such that $\varphi(Z) = \Lambda(Z)$ and
$A = (Y Z)^{-1}B\varphi(Y Z)$. From this equality follows that
\begin{equation}\label{f9}
A[A,Z^{-1}]C=ZA\varphi(Z)^{-1}=Y^{-1}B\varphi(Y),
\end{equation}
where $C = Z\Lambda(Z)^{-1} = Z\varphi(Z)^{-1}$ belongs to the center of ${\rm UT}_n(R)$. Looking to the image of equation (\ref{f9}) under the natural homomorphism $Z_{k+1}\to Z_{k+1}/Z_k$ we conclude that if
two elements $A,B \in Z_{k+1}~(k \neq0)$ are $\varphi$-conjugated, then their images under the natural homomorphism $Z_{k+1} \to Z_{k+1}/Z_k$ are $\varphi_k$-conjugated. Therefore for
all $k > 0$ we have $R(\varphi_k) < \infty$ and $R(\psi_k) < \infty$. By Proposition \ref{ind} we have $R(\psi_k) = \big|Z_{k+1}(G)/Z_k(G) : [e]_{\psi_k}\big|$, where $e$ is the unit element from $Z_{k+1}(G)/Z_k(G)$.
Since $R$ is a field of zero characteristic, $Z_{k+1}(G)/Z_k(G) = Z_{k+2}/Z_{k+1} = (R^+)^{k+2}$ is a
divisible abelian group and it has only one subgroup of finite index (itself). Therefore
for all $k > 0$ we have $R(\psi_k) = 1$ and from Corollary \ref{nilin} follows that $R(\psi) = 1$,
i. e. every element from $G$ is $\psi$-conjugated with the unit element, or equivalently,
every element from ${\rm UT}_n(R)$ is $\varphi$-conjugated with some element from the center of
${\rm UT}_n(R)$.\hfill$\square$

The following theorem follows from Corollary \ref{cor3} and Corollary \ref{cor4}.
\begin{theorem}Let $R$ be a field of zero characteristic and $n \geq 3$. Then the Reidemeister spectrum of ${\rm UT}_n(R)$ is $\{1,\infty\}$.
\end{theorem}
\noindent \textbf{Proof.} In the beginning of this section we understood that there exists a diagonal
automorphism $\varphi$ of ${\rm UT}_n(R)$ with $R(\varphi) = 1$. Also the group ${\rm UT}_n(R)$ has infinitely
many conjugacy classes, i.~e. $R(id) = \infty$. So, $\{1,\infty\}$ is a subset of $Spec_R({\rm UT}_n(R))$
and we need to prove that $1$ is the only finite number in $Spec_R({\rm UT}_n(R))$. 

Let $\varphi$ be an automorphism of ${\rm UT}_n(R)$ with $R(\varphi) < \infty$. By Proposition \ref{lev} there exist an inner automorphism $\varphi_A$, a diagonal automorphism $\psi_D$, a central
automorphism $\Lambda$, a ring automorphism $\Delta$ and an integer $m\in\{0,1\}$ such that $\varphi = \varphi_A\Lambda\sigma^m \psi_D\Delta$. By Proposition \ref{inn} we can assume that $\varphi_A = id$, so, $\varphi=\Lambda\sigma^m\psi_D\Delta$. Further it is a direct consequence of Corollary \ref{cor3} and Corollary \ref{cor4}.\hfill$\square$
\newpage

{\footnotesize
\begin{spacing}{0.5}

\end{spacing}}
\end{document}